\documentclass[11pt]{amsart}
\usepackage{amsmath,amsfonts,amssymb,amsxtra,latexsym,amscd,enumerate,amsthm}

\usepackage{latexsym}
\usepackage{epsfig}
\usepackage{verbatim}

\def\g{\gamma}
\def\G{\Gamma}
\def\a{\alpha}

\def\d{\delta}
\def\b{\beta}

\def\p{\Phi}

\def\ll{\Lambda}
\def\l{\lambda}
\def\o{\omega}

\def\n{\nu}
\def\R{\mathbb{R}}

\def\C{\mathbb{C}}
\def\Z{\mathbb{Z}}
\def\K{\mathcal{K}}

\def\p{\mathcal{P}}

\def\L{\mathcal{L}}
\def\V{\mathcal{V}}

\def\H{\mathcal{H}}

\def\N{\mathbb{N}}
\def\n{\mathcal{N}}
\def\f{\mathcal{F}}
\def\Z{\mathbb{Z}}
\def\I{\mathcal{I}}
\def\beq{\begin{equation}}
\def\eeq{\end{equation}}

\def\beq{\begin{equation}}
\def\eeq{\end{equation}}

\newtheorem{o0}{Observation}
\newtheorem{t1}{Theorem}
\newtheorem{l1}{Lemma}
\newtheorem{p1}{Proposition}

\begin{document}
\title[]{On the Boundedness of The Bilinear Hilbert Transform along ``non-flat" smooth curves.\\The Banach triangle case ($L^r,\: 1\leq r<\infty$).}

\author{Victor Lie}

\date{\today}

\address{Institute of Mathematics of the
Romanian Academy, Bucharest, RO 70700, P.O. Box 1-764, Romania.}

\address{Department of Mathematics, Purdue University, IN 47907 USA}

\email{vlie@purdue.edu}

\thanks{The author was supported by NSF grant DMS-1500958.}
\maketitle

\begin{abstract}
We show that the bilinear Hilbert transform $H_{\G}$ along curves $\Gamma=(t,-\gamma(t))$
with $\g\in\n\f^{C}$ is bounded from $L^{p}(\R)\times L^{q}(\R)\,\rightarrow\,L^{r}(\R)$ where $p,\,q,\,r$ are H\"older indices, i.e. $\frac{1}{p}+\frac{1}{q}=\frac{1}{r}$, with $1<p<\infty$, $1<q\leq\infty$ and $1\leq r<\infty$. Here $\n\f^{C}$ stands for a wide class of smooth ``non-flat" curves near zero and infinity whose precise definition is given in Section 2. This continues author's earlier work in \cite{L}, extending the boundedness range of $H_{\G}$ to any triple of indices $(\frac{1}{p},\,\frac{1}{q},\,\frac{1}{r'})$ within the Banach triangle. Our result is optimal up to end-points.
\end{abstract}
$\newline$
\section{\bf Introduction}

This paper, building upon the ideas in \cite{L}, continues the investigation of the boundedness properties of the bilinear Hilbert transform along curves. More precisely, if $\Gamma=(t, \, -\g(t))$ where here $\g$ is a suitable\footnote{In a sense that will be specified later.}
 smooth, non-flat curve near zero and infinity, we want to understand the behavior of the bilinear Hilbert transform along $\Gamma$ defined as
\beq\label{BHTc}
\eeq
$$H_{\G}: S(\R) \times S(\R)\longmapsto S'(\R)$$
$$H_{\G}(f,g)(x):= \textrm{p.v.}\int_{\R}f(x-t)g(x+\g(t))\frac{dt}{t}\:.$$
One can easily notice that taking $\g(t)=t$ we obtain the standard bilinear Hilbert transform. Thus, the problem considered in this paper is, in fact, the ``curved" analogue of the celebrated problem of providing $L^p$ bounds to the classical (``flat") bilinear Hilbert transform - the latter being solved in the seminal work of M. Lacey and C. Thiele (\cite{LT1}, \cite{LT2}).

It is worth noticing that similar studies regarding curved analogues of classical ``flat" objects arise naturally in harmonic analysis in various contexts. Indeed, recalling the discussion in the introduction of \cite{L}, a prominent such example  is given by the study of the boundedness of the linear Hilbert transform along curves given by
$$\H_{\G}:\: S(\R^n)\,\rightarrow\,S'(\R^n)\,,$$
$$\H_{\G}(f)(x):=p.v. \int_{\R} f(x-\G(t))\,\frac{dt}{t}\,,$$
where here $\Gamma:\:\R\,\rightarrow\,\R^n$, $n\geq 1$, is a suitable curve. This latter problem first appeared  in the work of Jones (\cite{Jo}) and  Fabes and Riviere (\cite{FR}) in connection with the analysis of the constant coefficient parabolic differential operators. The study of $H_{\G}$ was later extended to cover more general and diverse situations (\cite{NVWW} \cite{SW}, \cite{CNSW},\cite{C1}, \cite{C2}).

In the bilinear setting, the work on the ``curved" model referring to \eqref{BHTc} was initiated by X. Li in \cite{Li}. There, he showed that, for the particular case $\g(t)=t^d$ with $d\in\N$, $d\geq 2$, one has that
\beq\label{L2}
H_{\G}\:\:\textrm{maps}\:\:L^{2}(\R)\times L^{2}(\R)\,\rightarrow\,L^{1}(\R)\:.
\eeq
His proof relies on the concept of $\sigma-$uniformity, previously used in \cite{CLTT} and originating in the work of T. Gowers (\cite{G}).

In \cite{L}, the author proved that \eqref{L2} remains true for any curve $\g$ that belongs to $\n\f^{C}$ - a suitable  class\footnote{For its definition the reader is invited to consult Section 2. It is worth saying that $\n\f^{C}$ contains in particular the class of all polynomials without linear term - for more on this see Observation \ref{W}.} of smooth non-flat curves near zero and infinity. Our work greatly extended Li's result both qualitatively by significantly widening the class of curves, and quantitatively by revealing the scale type decay relative to the level sets of the multiplier's phase. Our proof is based on completely different methods, involving in a first instance a delicate analysis of the multiplier followed then by a special wave-packet discretization adapted to the two-directional oscillatory behavior of the phase. In the Appendix of the same paper, we also explained the main idea of how to upgrade the methods employed in \cite{Li} in order to be able to obtain the scale type decay.

Later, X. Li and L. Xiao, (\cite{LX}), relying heavily on \cite{L} in both the analysis treatment\footnote{See for example the decomposition of pp. 16 in \cite{LX} versus the one in Section 5.3. of \cite{L}, the treatment at pp. 29 in \cite{LX} versus the corresponding one at pp. 323 in \cite{L}, the perturbative strategy applied in analyzing the phase of the multiplier etc.} of the multiplier and the ``upgraded" $\sigma$-uniformity approach required for capturing the key necessary scale type decay revealed in \cite{L}, proved that if $P$ is any polynomial of degree $d\in\N$, $d\geq 2$ having no constant and no linear term \footnote{Any such polynomial $P$ is just a particular example of an element $\g\in\n\f$.} then taking $\g(t)=P(t)$ one has that $H_{\Gamma}$ maps boundedly $L^{p}(\R)\times L^{q}(\R)\:\:\mapsto\:\:L^{r}(\R)$ where $p,\,q,\,r$ obey $\frac{1}{p}+\frac{1}{q}=\frac{1}{r}$, $1<p<\infty$, $1<q<\infty$ and $\frac{d}{d-1}<r<\infty$. (The bounds here depend only on the degree $d$ but not on the coefficients of $P$). Finally, more recently, in \cite{GX}, the authors take the proof in \cite{L} and translate it in the $\sigma-$uniformity language used initially in \cite{Li}.

The present paper should be regarded as a natural continuation of \cite{L}. Relying on the author's previous methods, we extend the earlier results by showing that for any $\g\in\n\f$ the bounds on $H_{\Gamma}$ can be extended to cover the natural Banach triangle case. Our result is optimal up to end points. For a precise statement of our result as well as for a reminder of the definition of the class $\n\f$ of ``smooth non-flat" curves near zero and infinity one is invited to consult the next section.
$\newline$
\textbf{Acknowledgement.} I would like to thank my wife Anca for drawing Figure 1 in this paper.
$\newline$
\section{\bf Main results}
$\newline$

We start by recalling from \cite{L} the definition of the set $\n\f_{0}$ of all curves $\g$ which are smooth non-flat functions near the origin:

\begin{itemize}

\item \textit{smoothness, no critical points, variation (near origin)}
\beq\label{nocrtic0}
\eeq
$\exists\: \d>0$ (possibly depending on $\g$) and $\V_{\d}:=(-\d,\d)\setminus\{0\}$ such that $\g\in C^N(\V_{\d})$ $(N\geq4)$ and $|\g'|>0$ on $\V_{\d}$;
moreover
\beq\label{variation0}
\sup_{\a\in\R_{+}}\#\{j\in\Z_{+}\,|\,|2^{-j}\,\g'(2^{-j})|\in[\a,2\a]\}<\infty\:,
\eeq
where here $\Z_{+}:=\{j\in\Z\,|\,j\geq 0\}$.

\item \textit{asymptotic behavior (near origin)}

There exists $\{a_j\}_{j\in\N}\subset\R_{+}$ with $\lim_{j\rightarrow\infty}a_j=0$ such that:
$\newline$

For any $t\in I:=\{s\,|\,\frac{1}{4}\leq|s|\leq 4\}$ and $j\in \Z_{+}$ we have
\beq\label{asymptotic0}
\frac{\g(2^{-j}\,t)}{2^{-j}\,\g'(2^{-j})}=Q(t)\,+\,Q_j(t)\,,
\eeq
with $Q,\,Q_j\in C^N(I)$ and $\|Q_j\|_{C^N(I)}\leq a_j$.

For $s\in J=Q'(I)$ we require
\beq\label{asymptotic0gaminv}
\frac{(\g')^{-1}(s\,\g'(2^{-j}))}{(\g')^{-1}(\g'(2^{-j}))}=r(s)\,+\,r_j(s)\,,
\eeq
where $r,\,r_j\in C^{N-1}(J)$ with $\|r_j\|_{C^{N-1}(J)}\leq a_j\:.$

(The existence of $(\g')^{-1}$, the inverse of $\g'$, will be a consequence of the next hypothesis.)

\item \textit{non-flatness (near origin)}

The main terms in the asymptotic expansion obey
\beq\label{fstterm0}
 \inf_{t\in I} |Q''(t)|,\:\inf_{t\in J} |r'(t)|>c_{\g}>0\:,
\eeq
and
\beq\label{convdualphase0}
 \inf_{{t_1,\,t_2\in J}\atop{t_1\not=t_2}}\frac{|t_1\,r'(t_1)-t_2\,r'(t_2)|}{|t_1-t_2|}>c_{\g}\:.
\eeq
\end{itemize}

In a similar fashion one can define $\n\f_{\infty}$ - the class of smooth, non-flat near infinity functions $\g$ having the following properties:

\begin{itemize}

\item \textit{smoothness, no critical points, variation (near infinity)}
\beq\label{nocrticinf}
\eeq
$\exists\: \d>0$ (possibly depending on $\g$) and $\tilde{\V}_{\d}:=(-\infty,-\d)\cup (\d,\infty)$ such that $\g\in C^N(\tilde{\V}_{\d})$ $(N\geq4)$ and $|\g'|>0$ on $\tilde{\V}_{\d}$;
moreover
\beq\label{variationinft}
\sup_{\a\in\R_{+}}\#\{j\in\Z_{-}\,|\,|2^{-j}\,\g'(2^{-j})|\in[\a,2\a]\}<\infty\:,
\eeq
where here $\Z_{-}:=\{j\in\Z\,|\,j\leq 0\}$.

\item \textit{asymptotic behavior (near infinity)}

 There exists $\{\tilde{a}_j\}_{j\in \Z_{-}}\subset\R_{+}$ with $\lim_{j\rightarrow-\infty}\tilde{a}_j=0$ such that:

For any $t\in I:=\{s\,|\,\frac{1}{4}\leq|s|\leq 4\}$ and $j\in \Z_{-}$ we have
\beq\label{asymptoticinfty}
\frac{\g(2^{-j}\,t)}{2^{-j}\,\g'(2^{-j})}=\tilde{Q}(t)\,+\,\tilde{Q}_j(t)\,,
\eeq
with $\tilde{Q},\,\tilde{Q}_j\in C^N(I)$ and $\|\tilde{Q}_j\|_{C^N(I)}\leq \tilde{a}_j$.

For $s\in \tilde{J}=\tilde{Q}'(I)$ we require
\beq\label{asymptoticinftygaminv}
\frac{(\g')^{-1}(s\,\g'(2^{-j}))}{(\g')^{-1}(\g'(2^{-j}))}=\tilde{r}(s)\,+\,\tilde{r}_j(s)\,,
\eeq
where $\tilde{r},\,\tilde{r}_j\in C^{N-1}(\tilde{J})$ with $\|\tilde{r}_j\|_{C^{N-1}(\tilde{J})}\leq \tilde{a}_j\:.$

\item \textit{non-flatness (near infinity)}

The main terms in the asymptotic expansion obey
\beq\label{fstterminfty}
 \inf_{t\in I} |\tilde{Q}''(t)|,\:\inf_{t\in \tilde{J}} |\tilde{r}'(t)|>c_{\g}>0\:,
\eeq
and
\beq\label{convdualphaseinfty}
 \inf_{{t_1,\,t_2\in \tilde{J}}\atop{t_1\not=t_2}}\frac{|t_1\,\tilde{r}'(t_1)-t_2\,\tilde{r}'(t_2)|}{|t_1-t_2|}>c_{\g}\:.
\eeq
\end{itemize}

With this done, we set $$\n\f:=C(\R\setminus\{0\})\cap\n\f_{0}\cap\n\f_{\infty}$$ and $\n\f^{C}:=\n\f\,+\,Constant$.

\begin{o0}\label{W} Following \cite{L}, we list here some interesting features of the class $\n\f$:
\begin{itemize}
\item Any real polynomial of degree $\geq 2$ with no constant and no linear term belongs to $\n\f$;

\item More generally, any finite linear combination over $\R$ of terms of the from $t^{\a}$ with $\a\in (0,\infty)\setminus\{1\}$ belongs to $\n\f$;

\item Even more, any finite linear combination over $\R$ of terms of the form
$|t|^{\a}\,|\log |t||^{\b}$ with $\a,\,\b\in\R$ and $\a\notin\{-1,\,0,\,1\}$ is in $\n\f$;

\item If $\g\in\n\f$ then there exist $K_2\geq K_1>0$ and $C_2\geq C_1>0$ (all constants are allowed to depend on $\g$) such that for any $t\in V(0)\setminus\{0\}$ or $t\in V(\infty)$, respectively
\beq\label{growth}
\eeq
\begin{itemize}
\item either $K_2^{-1}\,|t|^{C_2}<|\g'(t)|<K_1^{-1}\,|t|^{C_1}$;
\item or $\frac{K_1}{|t|^{C_1}}<|\g'(t)|<\frac{K_2}{|t|^{C_2}}$.
\end{itemize}
$\newline$

Thus, if  $\g\in\n\f$, one has:
\beq\label{growth1}
\exists\:\:\lim_{{t\rightarrow 0}\atop{t\not=0}} |\g'(t)|\in\{0,\,\infty\}\:\:\:\textrm{and}\:\:\:\exists
\lim_{t\rightarrow \infty} |\g'(t)|\in\{0,\,\infty\}\:.
\eeq
\end{itemize}
\end{o0}

In \cite{L}, we proved the following result
$\newline$

\noindent\textbf{Theorem.} \textit{Let $\G=(t,-\gamma(t))$ be a curve such that $\g\in\n\f^{C}$. Recall the definition of the bilinear Hilbert transform $H_{\G}$ along the curve $\G$:
$$H_{\G}: S(\R) \times S(\R)\longmapsto S'(\R)$$
$$H_{\G}(f,g)(x):= \textrm{p.v.}\int_{\R}f(x-t)g(x+\g(t))\frac{dt}{t}\:.$$
Then $H_{\G}$ extends boundedly from $L^2(\R)\times L^2(\R)$ to $L^1(\R)$.}
$\newline$

In the present paper, we extend the boundedness range of the above theorem to the Banach triangle (see Figure 1):
$\newline$

\noindent\textbf{Main Theorem.} \textit{If $\g\in\n\f^{C}$ and $H_{\G}$ defined as above, we have that
\beq\label{bounds}
H_{\G}:\:\:\:\:\:L^{p}(\R)\times L^{q}(\R)\:\,\mapsto\:\, L^{r}(\R)\,,
\eeq
where the indices $p,\,q,\,r$ obey
\beq\label{holder}
\frac{1}{p}+\frac{1}{q}=\frac{1}{r}\,,
\eeq
with
\beq\label{ban}
1<p<\infty\,,\:\:\:1<q\leq\infty\,,\:\:\:\textrm{and}\:\:\:1\leq r<\infty\,.
\eeq}
$\newline$
\begin{o0}\label{optimal} The Banach range in our Main Theorem is optimal up to end-points. Indeed, let us define $\p_d$ be the class of all real polynomials of degree $d$ with no constant and no linear terms. Then, in \cite{Li} (pp. 9), it is shown that for any $d\in\N$, $d\geq 2$, there exists a polynomial $P_d\in \p_d$ such that for $\g=P_d$ one has that
\beq\label{cex}
\|H_{\G}\|_{L^{p}(\R)\times L^{q}(\R)\,\mapsto\,L^{r}(\R)}=\infty
\eeq
whenever $(p,\,q,\,r)$ obeys \eqref{holder} with $r<\frac{d-1}{d}$.

Since from Observation \ref{W} we deduce that for any $d\in\N$, $d\geq 2$, we have $\p_d\subset \n\f$ we conclude that in order for
\eqref{bounds} and \eqref{holder} to hold for \textit{any} $\g\in\n\f^{C}$ one must have $r\geq 1$.
\end{o0}

\begin{figure}[!h]
\begin{center}
\epsfig{file=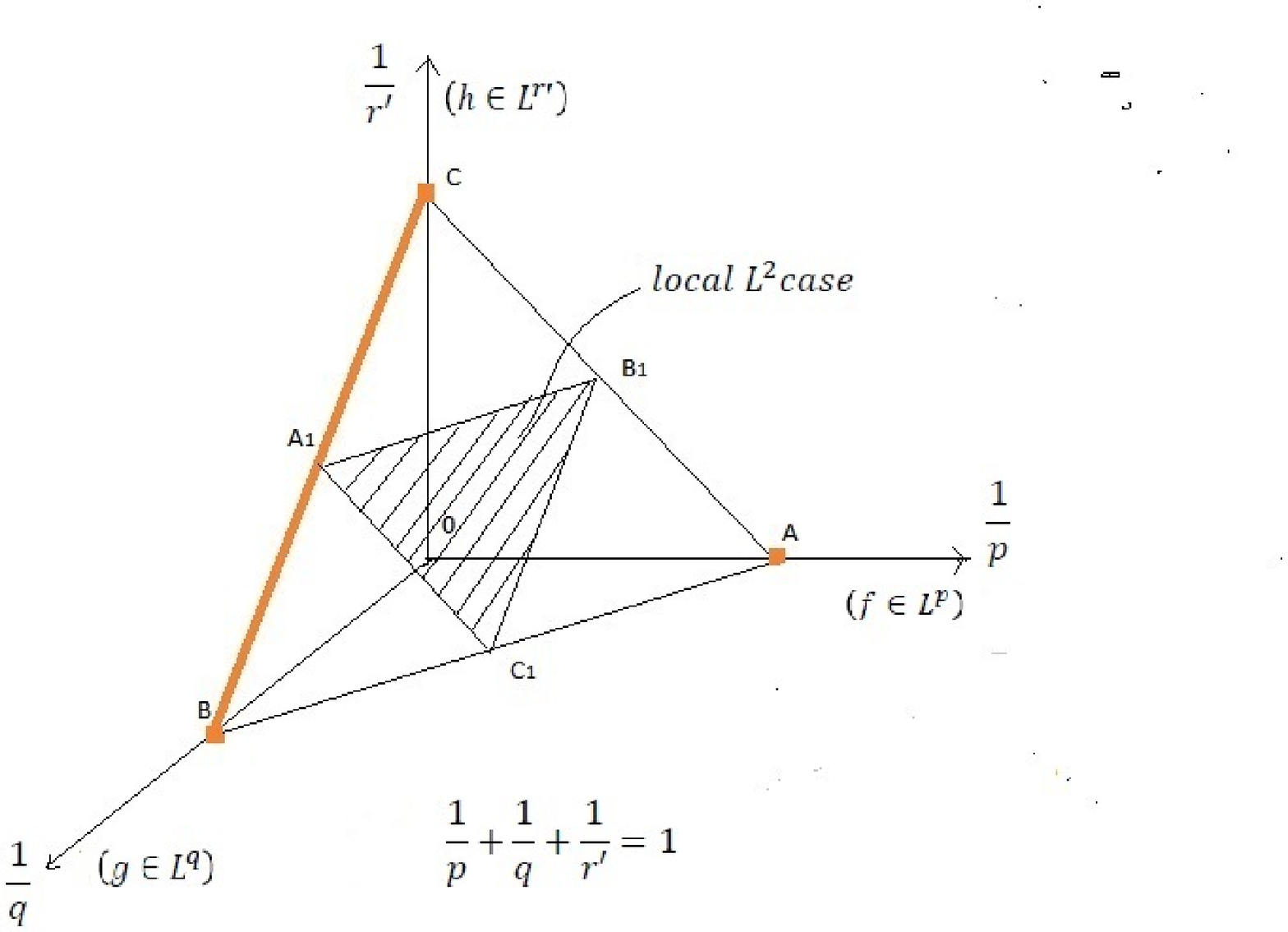,height=13cm}
\end{center}
\caption{\textbf{The boundedness range for the Bilinear Hilbert Transform $H_{\Gamma}$.}
In this figure we represent the bounds for our object viewed as a trilinear form defined by $\Lambda(f,g,h)=\int H_{\Gamma}(f,g)\,h$. Our Main Theorem states that $\Lambda$  maps boundedly $L^p\times L^q\times L^{r'}$ into $\C$ for all triples $(\frac{1}{p},\,\frac{1}{q},\,\frac{1}{r'})$ that belong to the region
$\Omega:=\textit{int}\,(\triangle ABC)\cup(AC)\cup(AB)\;.$
These bounds are optimal up to the bold boundaries defined by
$[BC]\cup\{A\}$.}
\end{figure}

\section{\bf Preparatives; isolating the main component of our operator}
$\newline$

In Section 3 of \cite{L}, after elaborated technicalities,  we proved that the study of our bilinear operator $H_{\G}(f,g)$, can be reduced to the corresponding study of the bilinear operator $T(f,g)$ defined as follows:
\beq\label{T}
T(f,g):=\sum_{j\in\Z}\sum_{m\in \N} T_{j,m}\,,
\eeq
with
\beq\label{pieces}
T_{j,m}(f,g)(x):= \int_{\R}\int_{\R}\hat{f}(\xi)\hat{g}(\eta)v_{j,m}(\xi,\eta) e^{i\xi x} e^{i\eta x} d\xi d\eta\:
 \eeq
and
\beq\label{vjm}
v_{j,m}(\xi,\eta):= 2^{-\frac{m}{2}}\,e^{i\,\varphi_{\xi,\eta}(t_c)}\,
\zeta\left(\frac{\eta\,\g'(2^{-j})}{2^{m+j}},\,\frac{\frac{\xi}{2^{m+j}}}{\frac{\eta\,\g'(2^{-j})}{2^{m+j}}}\right)\,
\phi\left(\frac{\xi}{2^{m+j}}\right)\,\phi\left(\frac{\eta\,\g'(2^{-j})}{2^{m+j}}\right)
 \eeq
where we have
\begin{itemize}
\item the function $\phi$ is smooth, compactly supported with
\beq\label{ssuport}
\textrm{supp}\,\phi\subset\{x\,|\,\frac{1}{10}<|x|<10\}\:.
\eeq

\item the phase of the multiplier is defined as
\beq\label{phase}
\varphi_{\xi,\eta}(t):=-\frac{\xi}{2^j}\,t+\eta\,\g(\frac{t}{2^j})\:.
 \eeq

\item for $\xi,\,\eta,\,j$ fixed (and based on the properties of $\g$) there exists exactly one critical point $t_c$
defined by
\beq\label{critpoint}
t_c=t_c(\xi\,,\eta,\,j)\in [2^{-k(\g)}, 2^{k(\g)}]\,\:\;\textrm{such that}\:\;\varphi'_{\xi,\eta}(t_c)=-\frac{\xi}{2^j}+\frac{\eta}{2^j}\,\g'(\frac{t_c}{2^j})=0\;,
\eeq
where here $k(\g)\in\N$ is an integer depending only on $\g$.

\item the function $\zeta$ obeys
\beq\label{zeta}
\zeta:\,[\frac{1}{10},10]\times[2^{-k(\g)},\,2^{k(\g)}]\,\rightarrow\,\R\:\:\:\textrm{with}\:\:\:\;\:\:\|\zeta\|_{C^{N-3}}\lesssim_{\g} 1\:.
\eeq
\end{itemize}

Now based on \eqref{growth1} and \eqref{growth}, wlog we can assume that

\beq\label{growth2}
\exists\:\:\lim_{{t\rightarrow 0}\atop{t\not=0}} \g'(t)=0\:\:\:\textrm{and}\:\:\:\exists
\lim_{t\rightarrow \infty} \g'(t)=\infty\:.
\eeq

The other cases can be treated similarly and we will not discuss them in detail here.

Setting
\beq\label{ph1}
\Psi_{\eta}(\xi)=-\varphi_{\xi,\eta}(t_c)\:,
\eeq
and following \cite{L}, we used the scaling symmetry in order to define the following operators:
\begin{itemize}

\item For $j>0$ (thus $2^{-j}\rightarrow 0$)

$$B_{j,m}(f(\cdot),g(\cdot))(x):=[\g'(2^{-j})]^{\frac{1}{2}}\,T_{j,m}\left( f(2^{m+j}\cdot), g(\frac{2^{m+j}}{[\g'(2^{-j})]}\cdot) \right)(\frac{x\,[\g'(2^{-j})]}{2^{m+j}})\:.$$
Remark that
\beq\label{poz0}
\eeq
$$B_{j,m}(f,g)(x)=2^{-\frac{m}{2}}\,[\g'(2^{-j})]^{\frac{1}{2}}\times$$
$$\int_{\R}\int_{\R}\hat{f}(\xi)\,\hat{g}(\eta)\,e^{i (\g'(2^{-j})\xi+\eta)\,x}\,e^{-i\,2^{m}\,2^j\,\Psi_{\frac{\eta}{\g'(2^{-j})}}(\xi)}\,
\zeta(\eta,\,\frac{\xi}{\eta})\,\phi(\xi)\,\phi(\eta)\,d\xi\,d\eta \:.$$

\item For $j<0$ (thus $2^{-j}\rightarrow \infty$)

$$B_{j,m}(f(\cdot),g(\cdot))(x):=
[\g'(2^{-j})]^{-\frac{1}{2}}\,T_{j,m}\left( f(2^{m+j}\cdot), g(\frac{2^{m+j}}{[\g'(2^{-j})]}\cdot) \right)(\frac{x}{2^{m+j}})\:.$$
As before, notice that
\beq\label{neg0}
\eeq
$$B_{j,m}(f,g)(x)=2^{-\frac{m}{2}}\,[\g'(2^{-j})]^{-\frac{1}{2}}\times$$
$$\int_{\R}\int_{\R}\hat{f}(\xi)\,\hat{g}(\eta)\,e^{i (\xi+\frac{\eta}{[\g'(2^{-j})]})\,x}\,e^{-i\,2^{m}\,2^j\,\Psi_{\frac{\eta}{\g'(2^{-j})}}(\xi)}\,
\zeta(\eta,\,\frac{\xi}{\eta})\,\phi(\xi)\,\phi(\eta)\,d\xi\,d\eta \:.$$
\end{itemize}
$\newline$

\begin{o0}\label{j}
In what follows we will focus on the case $j>0$ as the reasonings for the other case $j<0$ can be treated in a similar fashion.
\end{o0}
$\newline$
$\newline$
\textbf{Main reductions}
$\newline$

For $j\in\N$, using \eqref{growth2} and $\g\in\n\f$, one follows the reasonings from Section 5 in \cite{L} and successively simplifies the structure of $B_{j,m}$ as follows:
\begin{itemize}
\item for $s\in J$, recalling the definition of $r$ and $r_j$ in \eqref{asymptotic0gaminv}, we set
\beq\label{Psi}
R(s):=\int_{1}^{s} r(u)\, du\:\:\:\textrm{and}\:\:\:R_j(s):=\int_{1}^{s} r_j(u)\, du\,,
\eeq
and first notice that
\beq\label{Psi}
\begin{array}{rl}
&2^j\,\Psi_{\frac{\eta}{\g'(2^{-j})}}(\xi)=
 \eta\,R(\frac{\xi}{\eta})\,+\,\eta\,R_j(\frac{\xi}{\eta})\,.
\end{array}
\eeq
Based on the properties of $\g$ and on the assumptions on $\zeta$ it is enough to treat the term
\beq\label{Bjmintdec}
\eeq
$$B_{j,m}(f,g)(x)\approx
 2^{-\frac{m}{2}}[\g'(2^{-j})]^{\frac{1}{2}} $$
 $$\int_{\R}\int_{\R}\hat{f}(\xi)\,\hat{g}(\eta)\,e^{i (\g'(2^{-j})\xi+\eta)\,x}\,
e^{-i\,2^m\,\eta\,R(\frac{\xi}{\eta})}\,\phi(\xi)\,\phi(\eta)\,d\xi\,d\eta \:.$$
and thus
\beq\label{Tjm}
\eeq
$$T_{j,m}(f,g)(x)\approx$$
$$2^{-\frac{m}{2}}
 \int_{\R}\int_{\R}\hat{f}(\xi)\,\hat{g}(\eta)\,e^{i (\xi+\eta)\,x}\,
e^{-i\,\frac{\g'(2^{-j})}{2^{j}}\,\eta\,R(\frac{\xi}{\g'(2^{-j})\,\eta})}\,\phi(\frac{\xi}{2^{m+j}})\,
\phi(\frac{\g'(2^{-j})\,\eta}{2^{m+j}})\,d\xi\,d\eta \:.$$

\item Because the phase of the multiplier has roughly the size $2^m$ with the dominant factor in the variable $\eta$, we further divide the $\eta-$support in $\approx 2^{m}$ intervals of equal size. Thus, after running the same approximation algorithm as in \cite{L}, we deduce that the main component of $T_{j,m}$ is given by
\beq\label{mainTjm}
\eeq
$$\tilde{T}_{j,m}(f,g)(x):=$$
$$2^{-\frac{m}{2}}
\sum_{p_0=2^m}^{2^{m+1}} \int_{\R}\int_{\R}\hat{f}(\xi)\,\hat{g}(\eta)\,e^{i (\xi+\eta)\,x}\,
e^{-i\,p_0\,R(\frac{\xi}{2^j\,p_0})}\,\phi(\frac{\xi}{2^{m+j}})\,
\phi(\frac{\g'(2^{-j})\,\eta}{2^{j}}-p_0)\,d\xi\,d\eta \:.$$
\end{itemize}
$\newline$

For notational simplicity, we will allow a small abuse, and redenote $\tilde{T}_{j,m}$ as just simply $T_{j,m}$.
$\newline$

Let us introduce now several notations:
\beq\label{a}
\phi_{k}(\xi):=\phi(\frac{\xi}{2^{k}})\,,\:\:\:\:\:k\in\Z\,;
\eeq
\beq\label{b}
\phi_{j,p_0}(\eta):=\phi(\frac{\g'(2^{-j})\,\eta}{2^{j}}-p_0)\,,\:\:\:\:\:j\in\N,\:p_0\in [2^m,\,2^{m+1})\cap\N\,;
\eeq
\beq\label{c}
\psi_{m,p_0,j}(\xi):=2^{-\frac{m}{2}}\,e^{-i\,p_0\,R(\frac{\xi}{2^j\,p_0})}\,\phi(\frac{\xi}{2^{m+j}}),\:\:\:\:\:m\in\N,\:p_0\in [2^m,\,2^{m+1})\cap\N,\:j\in\N\,.
\eeq
It will be useful to record the following identity
\beq\label{ftr}
\check{\psi}_{m,p_0,j}(x)=2^{j}\,e^{i\,\vartheta_{p_0}(2^{j}\,x)}\,\phi^{*}(\frac{p_0}{2^m}\,r^{-1}(2^j\,x))\,+\,\textrm{Error Term}\:,
\eeq
where here
\beq\label{phase}
\vartheta_{p_0}(2^{j}\,x):=p_0\,\int_{0}^{r^{-1}(2^{j}\,x)} t\,r'(t)\,dt\,,
\eeq
$\phi^{*}\in C_0^{\infty}(\R)$ is a smooth function supported away from the origin with $\|\phi^{*}\|_{C^{100}}\lesssim 1$ and finally, the Error represents a smooth, fast decaying term relative to the $m$-parameter.

Also, throughout the paper, we have $p,\,q>1$, $\frac{1}{p}+\frac{1}{q}=\frac{1}{r}$ with $1\leq r <\infty$ and $\frac{1}{r'}=1-\frac{1}{r}$.

Now for $f\in L^{p}(\R)$, $g\in L^{q}(\R)$ and $h\in L^{r'}(\R)$ (recall $j\in\N$) we set
\beq\label{ljm}
\ll_{j,m}(f,g,h):=\int_{\R} T_{j,m}(f,g)(x)\,h(x)\,dx\:.
\eeq
From the above considerations, using relation \eqref{growth2} and Observation \ref{j}, we deduce that\footnote{Throughout this paper we will ignore the contribution of the error term in \eqref{ftr} as its treatment was carefully presented in \cite{L}. Also, for notational convenience, we make a small abuse by letting the same function $\phi_{j,p_0}$ represent the Fourier localization of both $g$ and $h$; in reality the support of $\hat{h}$ should be slightly (e.g. twice) larger than the support of $\hat{g}$.}
\beq\label{ljmexp}
\ll_{j,m}(f,g,h)=\sum_{p_0=2^m}^{2^{m+1}}\int_{\R} (f*\check{\phi}_{m+j}*\check{\psi}_{m,p_0,j})(x)\,(g*\check{\phi}_{j,p_0})(x)\,(h*\check{\phi}_{j,p_0})(x)\,dx\:.
\eeq
We now define
\beq\label{Lambdamp}
\ll_{m}^{+}(f,g,h):=\sum_{j\in \N}\ll_{j,m}(f,g,h)\,.
\eeq
A similar decomposition can be done for treating the case $j<0$. In this situation the analogue of \eqref{Lambdamp} will read
\beq\label{Lambdamp1}
\ll_{m}^{-}(f,g,h):=\sum_{j\in \Z\setminus\N}\ll_{j,m}(f,g,h)\,,
\eeq
where for this case $\ll_{j,m}$ is an appropriate adaptation of \eqref{ljmexp} to the context of \eqref{neg0} instead of \eqref{poz0}. For the brevity of the exposition, we won't insist more on this.

With all these done, our main focus will be on understanding the main properties of the trilinear form
\beq\label{Lambda}
\ll(f,g,h):=\sum_{m\in\N}\ll_{m}(f,g,h)\,,
\eeq
with
\beq\label{Lambdam}
\ll_{m}(f,g,h):=\ll_{m}^{+}(f,g,h)\,+\,\ll_{m}^{-}(f,g,h)\:.
\eeq

Based on the reductions presented above, we deduce that for proving our Main Theorem it is enough to show
$\newline$

\begin{t1}\label{Ll1} The form $\ll$ initially defined on $S(\R) \times S(\R)\times S(\R)\longmapsto \C$ by \eqref{Lambda} obeys the bounds:
\beq\label{bound}
\ll:\:L^{p}(\R)\times L^{q}(\R)\times L^{r'}(\R)\:\:\rightarrow\:\:\C
\eeq
 where the indices $p,\,q,\,r'$ satisfy
$$\frac{1}{p}+\frac{1}{q}+\frac{1}{r'}=1\,,$$ with
$$1<p<\infty\,,\:\:\:1<q\leq\infty\,,\:\:\:\textrm{and}\:\:\:1\leq r<\infty\,.$$
\end{t1}
$\newline$

Give relation \eqref{Lambda}, we further notice that Theorem \ref{Ll1} follows from the theorem below\footnote{For a helpful geometric perspective, the reader is invited to consult Figure 1.} after applying real interpolation methods and a telescoping sum argument (see Section 5):
$\newline$

\begin{t1}\label{estimabac} \textit{Let $1<p<\infty$. Then the following estimates hold \footnote{Throughout the paper $p^*=\min\{p,\,p'\}\;.$}}:
\beq\label{ac}
\underline{\textrm{Edge}\:\:(AC)}\,:\:\:\:\:\:\:\:\:\:\:\:|\ll_{m}(f,g,h)|\lesssim_{p} (1+m^{1-\frac{2}{p}})\,2^{-\frac{m}{10}(1-\frac{1}{p^{*}})}\,\|f\|_{p}\,\|g\|_{\infty}\,\|h\|_{p'}\;,
\eeq
and
\beq\label{ab}
\underline{\textrm{Edge}\:\:(AB)}\,:\:\:\:\:\:\:\:\:\:\:\:|\ll_{m}(f,g,h)|\lesssim_{p} (1+m^{1-\frac{2}{p}})\,2^{-\frac{m}{10}(1-\frac{1}{p^{*}})}\,\|f\|_{p}\,\|g\|_{p'}\,\|h\|_{\infty}\;.
\eeq
\end{t1}
$\newline$

Finally, based on Observation \ref{j}, it will be enough to prove Theorem \ref{estimabac} above for the ``positive half-line" forms $\ll_{m}^{+}(f,g,h)$ with $m\in\N$.
$\newline$

\section{Controlled bounds for $\ll_{m}^{+}(f,g,h)$ on the edge $(AC)$.}
$\newline$

 In this section, we focus on providing ``slowly increasing" bounds\footnote{Relative to the $m-$parameter.} for our form $\ll_{m}^{+}(f,g,h)$ on the edge $(AC)$ (see Figure 1). More precisely, our goal is to get the following tamer result:
$\newline$

\begin{t1}\label{tame} Let $1<p<\infty$. Then the following estimate holds:
\beq\label{ac}
\underline{\textrm{Edge}\:\:(AC)}\,:\:\:\:\:\:\:\:\:\:\:\:|\ll_{m}^{+}(f,g,h)|\lesssim_p (1+m^{\frac{2}{p'}-1})\,\|f\|_{p}\,\|g\|_{\infty}\,\|h\|_{p'}\;.
\eeq
\end{t1}
$\newline$

\subsection{Estimates for the edge $(AB_1]$.} The main result of this subsection is

\begin{p1}\label{estimab1} If $1<p\leq 2$, the following holds:
\beq\label{pointAB1}
\underline{\textrm{Edge}\:\:(AB_1]}\,:\:\:\:\:\:\:\:\:\:\:\:|\ll_{m}^{+}(f,g,h)|\lesssim_{p}  \|f\|_{p}\, \|g\|_{\infty}\,\|h\|_{p'}\,.
\eeq
\end{p1}
$\newline$

\begin{proof}

As a consequence of \eqref{ljmexp}, we deduce that
\beq\label{shapelambda1}
\eeq
$$|\ll_{m}^{+}(f,g,h)|\lesssim\sum_{j\in\N}\int_{\R} |(f*\check{\phi}_{m+j}(y)|$$
$$\times\,\int_{\R}\left(\sum_{p_0=2^m}^{2^{m+1}}|(g*\check{\phi}_{j,p_0})(x)\,(h*\check{\phi}_{j,p_0})(x)|\right)
\,2^{j}\,|\phi^{*}(r^{-1}(2^j\,(x-y))|\,dx\,dy\:.$$

Now, letting $M$ stand for the standard Hardy-Littlewood maximal function, we deduce the key relation
\beq\label{key1}
|\ll_{m}^{+}(f,g,h)|\lesssim
\sum_{j\in\N}\int_{\R} |(f*\check{\phi}_{m+j}(y)|\:\:
M\left(\sum_{p_0=2^m}^{2^{m+1}}|(g*\check{\phi}_{j,p_0})\,(h*\check{\phi}_{j,p_0})|\right)(y)\,dy\:.
\eeq

Now, from \eqref{key1}, we deduce that
$$|\ll_{m}^{+}(f,g,h)|\lesssim$$
$$\int_{\R}\left( \sum_{j\in\N}|(f*\check{\phi}_{m+j}(y)|^2\right)^{\frac{1}{2}}\:
\left( \sum_{j\in\N} \left|M\left(\sum_{p_0=2^m}^{2^{m+1}}|(g*\check{\phi}_{j,p_0})\,(h*\check{\phi}_{j,p_0})|\right)(y)\right|^2\right)^{\frac{1}{2}}\,dy$$
$$\lesssim \left\| \left( \sum_{j\in\N}|(f*\check{\phi}_{m+j}(y)|^2\right)^{\frac{1}{2}}\right\|_{p}\,
\left\|\left( \sum_{j\in\N} \left|M\left(\sum_{p_0=2^m}^{2^{m+1}}|(g*\check{\phi}_{j,p_0})\,(h*\check{\phi}_{j,p_0})|\right)(y)\right|^2\right)^{\frac{1}{2}}\right\|_{p'}$$
$$\lesssim \|f\|_{p}\,\left\|\left( \sum_{j\in\N}\left(\sum_{p_0=2^m}^{2^{m+1}}|(g*\check{\phi}_{j,p_0})\,(h*\check{\phi}_{j,p_0})|\right)^2\right)^{\frac{1}{2}}\right\|_{p'}\:,$$
where for the last relation we used standard Littelwood-Paley theory (for providing bounds on the square function for $f$) and Fefferman-Stein's inequality (\cite{FS}) (for the term involving the functions $g$ and $h$).

Now we make use of the following observation:
\beq\label{L2infcont}
\sum_{p_0=2^m}^{2^{m+1}}|(g*\check{\phi}_{j,p_0})|^2\lesssim \|g\|_{\infty}^2\,.
\eeq

Indeed, this is a simple consequence of the fact that
\beq\label{conv}
(g*\check{\phi}_{j,p_0})(x)=\int_{\R} g(x-\frac{\g'(2^{-j})}{2^j}\,y)\,\check{\phi}(y)\,e^{i\,p_0\,y}\,dy\:.
\eeq

Thus applying Cauchy-Schwarz inequality we get
\beq\label{CS}
\left(\sum_{p_0=2^m}^{2^{m+1}}|(g*\check{\phi}_{j,p_0})\,(h*\check{\phi}_{j,p_0})|\right)^2\lesssim \|g\|_{\infty}^2\,\sum_{p_0=2^m}^{2^{m+1}}|(h*\check{\phi}_{j,p_0})|^2\,.
\eeq
Inserting \eqref{CS} in the last estimate on $|\ll_{m}^{+}(f,g,h)|$ we get:
\beq\label{L}
|\ll_{m}^{+}(f,g,h)|\lesssim  \|f\|_{p}\, \|g\|_{\infty}\,\left\|\left( \sum_{j\in\N}\sum_{p_0=2^m}^{2^{m+1}}|h*\check{\phi}_{j,p_0}|^2\right)^{\frac{1}{2}}\right\|_{p'}\:.
\eeq

Define now
$$A_{j,p_0}:=\left[\frac{p_0-10}{2^{-j}\,\g'(2^{-j})},\,\frac{p_0-\frac{1}{10}}{2^{-j}\,\g'(2^{-j})}\right]
\cup \left[\frac{p_0+\frac{1}{10}}{2^{-j}\,\g'(2^{-j})},\,\frac{p_0+10}{2^{-j}\,\g'(2^{-j})}\right]\,.$$
From \eqref{ssuport} and \eqref{b} we notice that
$$\textrm{supp} \:\phi_{j,p_0}\subseteq A_{j,p_0}\:.$$

Using now the key assumption \eqref{variation0}, we deduce that $\{A_{j,p_0}\}_{j,p_0}$ have the finite intersection property, that is
\beq\label{finint}
\sum_{j\in\N}\sum_{p_0=2^m}^{2^{m+1}}\chi_{A_{j,p_0}}\lesssim_{\g} 1\,.
\eeq
Now, due to \eqref{finint} and the fact that $p'\geq2$, we are precisely in the setting of Rubio de Francia's inequality (\cite{RF}):
\beq\label{RF}
\left\|\left( \sum_{j\in\N}\sum_{p_0=2^m}^{2^{m+1}}|h*\check{\phi}_{j,p_0}|^2\right)^{\frac{1}{2}}\right\|_{p'}\lesssim_{p'} \|h\|_{p'}\:.
\eeq
Thus, putting together \eqref{CS}, \eqref{L}, \eqref{RF} we conclude that \eqref{pointAB1} holds.
\end{proof}

$\newline$

\subsection{Estimates for the edge $[B_1 C)$.} In this subsection we will focus on proving

\begin{p1}\label{estimb1c} If $p\geq 2$, the following holds:
\beq\label{pointB1C}
\underline{\textrm{Edge}\:\:[B_1 C)}\,:\:\:\:\:\:\:\:\:\:\:\:|\ll_{m}^{+}(f,g,h)|\lesssim_{p} m^{\frac{2}{p'}-1}\,\|f\|_{p}\,\|g\|_{\infty}\,\|h\|_{p'}\;.
\eeq
\end{p1}
$\newline$

The proof of Proposition \ref{estimb1c} will be decomposed in two main steps encoded in the two propositions below:

\begin{p1}\label{Contf} Let $p\geq2$. Then, with the previous notations, we have
\beq\label{Cf}
\left\|\,\left(\sum_{j\in\N}\sum_{p_0=2^m}^{2^{m+1}}|f*\check{\phi}_{m+j}*\check{\psi}_{m,p_0,j}(x)|^2\right)^{\frac{1}{2}}\right\|_{p}\lesssim_{p} m^{\frac{2}{p'}-1}\, \|f\|_{p}\;.
\eeq
\end{p1}

\begin{p1}\label{Contgh} If $1<p<\infty$ and $\frac{1}{p}+\frac{1}{p'}=1$ then the following holds:
\beq\label{Cgh}
\left\|\left(\sum_{j\in\N}\sum_{p_0=2^m}^{2^{m+1}}|(g*\check{\phi}_{j,p_0})(x)\,(h*\check{\phi}_{j,p_0})(x)|^{2}\right)^{\frac{1}{2}}\right\|_{p'}
\lesssim_{p'}\|g\|_{\infty}\,\|h\|_{p'}\;.
\eeq
\end{p1}

Assuming for the moment that \eqref{Cf} and \eqref{Cgh} hold we pass to the
$\newline$

\noindent\textbf{Proof of Proposition \ref{estimb1c}.}
$\newline$

By Cauchy-Schwarz we have
$$|\ll_{m}^{+}(f,g,h)|\leq$$
$$\sum_{j\in\N}\int_{\R}\left(\sum_{p_0=2^m}^{2^{m+1}}|(g*\check{\phi}_{j,p_0})(x)\,(h*\check{\phi}_{j,p_0})(x)|^{2}\right)^{\frac{1}{2}}
\,\left(\sum_{p_0=2^m}^{2^{m+1}}|f*\check{\phi}_{m+j}*\check{\psi}_{m,p_0,j}(x)|^2\right)^{\frac{1}{2}}\,dx\:,$$
which after a second application of Cauchy-Schwarz and then H\"older inequality becomes:
\beq\label{key2}
\eeq
$$|\ll_{m}^{+}(f,g,h)|\leq$$
$$\left\|\left(\sum_{j\in\N}\sum_{p_0=2^m}^{2^{m+1}}|(g*\check{\phi}_{j,p_0})(x)\,(h*\check{\phi}_{j,p_0})(x)|^{2}\right)^{\frac{1}{2}}\right\|_{p'}
\left\|\,\left(\sum_{j\in\N}\sum_{p_0=2^m}^{2^{m+1}}|f*\check{\phi}_{m+j}*\check{\psi}_{m,p_0,j}(x)|^2\right)^{\frac{1}{2}}\right\|_{p}\:.$$

Now applying Propositions \ref{Contf} and \ref{Contgh} we conclude that \eqref{pointB1C} holds.

We pass now to the
$\newline$

\noindent\textbf{Proof of Proposition \ref{Contf}.}
$\newline$

This proof will be decomposed in three main components.

The first one relies on the following key aspect: the
\textit{cancelation} offered by the phase of $\check{\psi}_{m,p_0,j}$; indeed, we have that
$\newline$

\begin{l1}\label{canc} The following holds:
\beq\label{claimc}
\sum_{p_0=2^m}^{2^{m+1}}|f*\check{\phi}_{m+j}*\check{\psi}_{m,p_0,j}(x)|^2\lesssim \int_{\R}|(f*\check{\phi}_{m+j})(y)|^2\, 2^j\,\nu(2^j\,(x-y))\,dy\:,
\eeq
where here $\nu\in C_0^{\infty}(\R)$ with $0\not\in\textrm{supp}\,\nu$ and $\|\nu\|_{C^{10}(\R)}\lesssim 1$.
\end{l1}

\begin{proof}

For fixed $m\in \N,\,j\in\N$ set now $w(y):=(f*\check{\phi}_{m+j})(\frac{y}{2^j})$. Then \eqref{claimc} turns into
\beq\label{claim1}
\sum_{p_0=2^m}^{2^{m+1}}|\int_{\R} w_x(y)\,e^{-i\,p_0\,\int_{0}^{r^{-1}(y)} t\,r'(t)\,dt}\,\phi^{*}(\frac{p_0}{2^m}\,r^{-1}(y))\,dy|^2\lesssim \int_{\R}|w_x(y)|^2\, \nu(y)\,dy\:.
\eeq
where here $x\in\R$ is considered fixed and we set $w_x(y):=w(2^j\,x-y)$.

Recall that from \eqref{phase}, we have $p_0\,\int_{0}^{r^{-1}(y)} t\,r'(t)\,dt=p_0\,\vartheta_{1}(y):=p_0\,\vartheta(y)$
Also recall that for expression \eqref{claim1} to make sense one should also have\footnote{The constant $C_\g$ are allowed to
change from line to line.}
\beq\label{cond}
\frac{1}{C_{\g}}<|y|<C_{\g}\:\:\:\textrm{for some fixed}\:C_{\g}>1\;.
\eeq
Set then $\mu\in C_{0}^{\infty}(\R)$ with $\textrm{supp}\,\mu\subseteq[\frac{1}{2 C_{\g}}, 2C_{\g}]$ and $\mu|_{[\frac{1}{ C_{\g}}, C_{\g}]}=1$.

Now in order to prove \eqref{claim1} it is enough to show the dual statement
\beq\label{claim2}
\int_{\R}|\sum_{p_0=2^m}^{2^{m+1}}\,a_{p_0}\,e^{i\,p_0\,\vartheta(y)}\,\phi^{*}(\frac{p_0}{2^m}\,r^{-1}(y))\,\mu (y)\,|^2\,dy\lesssim
\sum_{p_0=2^m}^{2^{m+1}} |a_{p_0}|^2 \:.
\eeq

For this we analyze the generic interaction
\beq\label{gen}
E_{p_0,q_0}:=\int_{\R}e^{i\,(p_0-q_0)\,\vartheta(y)}\,\phi^{*}(\frac{p_0}{2^m}\,r^{-1}(y))\,\overline{\phi^{*}}(\frac{q_0}{2^m}\,r^{-1}(y))\,\mu (y)\,\mu(y)\,dy\:.
\eeq

Now since $\vartheta'(y)=r^{-1}(y)$ has obviously no roots inside the support of $\mu$ we deduce that for $2^{m-10}<p_0, q_0< 2^{m+10}$ one has
\beq\label{gen1}
|E_{p_0,q_0}|\lesssim \frac{1}{1+|p_0-q_0|^2}\:.
\eeq
Thus
\beq\label{CSc}
\eeq
$$\int_{\R}|\sum_{p_0=2^m}^{2^{m+1}}\,a_{p_0}\,e^{i\,p_0\,\vartheta(y)}\,\phi^{*}(\frac{p_0}{2^m}\,r^{-1}(y))\,\mu (y)\,|^2\,dy$$
$$\lesssim\sum_{p_0,q_0=2^m}^{2^{m+1}} |a_{p_0}|\,|a_{q_0}|\,\frac{1}{1+|p_0-q_0|^2}\lesssim \sum_{p_0=2^m}^{2^{m+1}} |a_{p_0}|^2\:,$$
proving this way our claim \eqref{claim2} and hence \eqref{claimc}.
\end{proof}

We pass now to the second ingredient:
$\newline$

\begin{l1}\label{max}
Let $u\in L^1_{loc}(\R)$ and $\nu$ as before. Then, preserving the previous notations, we have:
\beq\label{claim4}
\begin{array}{cl}
I_{m,j}(u)(x):=\int|(u*\check{\phi}_{j+m})(y)|^2\,2^{j}\,\nu(2^j(x-y))\,dy\\\\
\lesssim \frac{1}{2^m}\,\sum_{l=2^{m-1}}^{2^{m+1}}\,|M(u)(x-\frac{l}{2^{m+j}})|^2\:.
\end{array}
\eeq
\end{l1}

\begin{proof}

 From the assumptions made in our hypothesis, using a smooth, compactly supported partition of unity, we have
\beq\label{dec}
\nu(y)=\sum_{l=2^{m-1}}^{2^{m+1}} \psi(2^m y-l)\:.
\eeq
From this, we deduce
\beq\label{desc}
\begin{array}{cl}
I_{m,j}(u)(x)=\int_{\R} \left(\int_{\R} u(s)\,2^{m+j}\,\check{\phi}(2^{m+j}(y-s))\,ds\right)\\\\
\times\left(\int_{\R}\overline{u(t)}\,2^{m+j}\,\overline{\check{\phi}(2^{m+j}(y-t))}\,dt\right)
\,\left(2^{j}\,\sum_{l=2^{m-1}}^{2^{m+1}} \psi(2^{m+j}\,(x-y)-l)\right)\,dy\\\\
=\frac{1}{2^m}\,\int_{\R\times\R} u(s)\,\overline{u(t)}\,\K_{m,j}(s,t)\,ds\,dt\,,
\end{array}
\eeq
where we set
\beq\label{defK}
\K_{m,j}(s,t):=2^{3(m+j)}\,\sum_{l=2^{m-1}}^{2^{m+1}}
\int_{\R}\check{\phi}(2^{m+j}(y-s))\,\overline{\check{\phi}(2^{m+j}(y-t))}\,\psi(2^{m+j}\,(x-y)-l)\,dy
\eeq
From the smooth behavior of $\check{\phi}$ we deduce
\beq\label{boundK}
|\K_{m,j}(s,t)|\lesssim 2^{2(m+j)}\,\sum_{l=2^{m-1}}^{2^{m+1}}
|\check{\phi}|(2^{m+j}(x-s)-l)\,|\check{\phi}|(2^{m+j}(x-t)-l)\,.
\eeq
Thus, from \eqref{desc} and \eqref{boundK} we have
\beq\label{cont}
|I_{m,j}(x)|\lesssim \frac{1}{2^m}\,\sum_{l=2^{m-1}}^{2^{m+1}}\,\left(\int_{\R} |u(s)|\,2^{m+j}\,
|\check{\phi}|(2^{m+j}(x-s)-l)\,ds\,\right)^2\,,
\eeq
from which we conclude
\beq\label{conc}
|I_{m,j}(x)|\lesssim \frac{1}{2^m}\,\sum_{l=2^{m-1}}^{2^{m+1}}\,[M (|u(\frac{\cdot}{2^{m+j}})|)(2^{m+j}x-l)]^2\,,
\eeq
which proves \eqref{claim4}.

\end{proof}

Combining now Lemmas \ref{canc} and \ref{max}, we deduce that
\beq\label{prop1}
\eeq
$$\left\|\,\left(\sum_{j\in\N}\sum_{p_0=2^m}^{2^{m+1}}\left|f*\check{\phi}_{m+j}*\check{\psi}_{m,p_0,j}(x)\right|^2\right)^{\frac{1}{2}}\right\|_{p}$$
$$\lesssim \left\|\,\left(\sum_{j\in\N}\frac{1}{2^m}\,\sum_{l\approx 2^m}\,\left|M(f*\check{\phi}_{m+j})(\cdot-\frac{l}{2^{m+j}})\right|^2\right)^{\frac{1}{2}}\right\|_{p}\:.$$
On the other hand, from Fefferman-Stein's inequality, (\cite{FS}), we have
\beq\label{FS}
\eeq
$$\left\|\,\left(\sum_{j\in\N}\frac{1}{2^m}\,\sum_{l\approx 2^m}\,\left|M(f*\check{\phi}_{m+j})(\cdot-\frac{l}{2^{m+j}})\right|^2\right)^{\frac{1}{2}}\right\|_{p}$$
$$\lesssim \left\|\,\left(\sum_{j\in\N}\frac{1}{2^m}\,\sum_{l\approx2^m}\,\left|(f*\check{\phi}_{m+j})(\cdot-\frac{l}{2^{m+j}})\right|^2\right)^{\frac{1}{2}}\right\|_{p}\;.$$

However, applying Fubini and Jensen's inequality (at this last point we use the key fact that $p\geq 2$)
$$\left\|\,\left(\sum_{j\in\N}\frac{1}{2^m}\,\sum_{l\approx2^m}\,\left|(f*\check{\phi}_{m+j})(\cdot-\frac{l}{2^{m+j}})\right|^2\right)^{\frac{1}{2}}\right\|_{p}^{p}$$
$$\lesssim\frac{1}{2^m}\,\sum_{l\approx2^m}\left\|\,\left(\sum_{j\in\Z}\,\left|(f*\check{\phi}_{m+j})(\cdot-\frac{l}{2^{m+j}})\right|^2\right)^{\frac{1}{2}}\right\|_{p}^{p}$$
Now, we appeal to the last ingredient in order to prove Proposition \ref{Contf}:
$\newline$

\begin{l1}\label{sqf}
Take $l\in\Z$ and define the $l-$shifted square function
\beq\label{Sqf}
S_{l}f(x):=\left(\sum_{j\in\Z}\left|(f*\check{\phi}_{j})(x-\frac{l}{2^{j}})\right|^2\right)^{\frac{1}{2}}\;.
\eeq
 Then, for any $1<q<\infty$, one has
\beq\label{Sqfp}
\|S_l f\|_q\lesssim_q (\log (|l|+10))^{\frac{2}{q^*}-1}\,\|f\|_q\;.
\eeq
\end{l1}
$\newline$

\begin{proof}

This statement relies on interpolation and Calderon-Zygmund theory adding to the classical theory of the $L^p$ boundedness of the (classical) square function a small extra twist in order to deal with the translation symmetry encoded in the $l$-parameter. In a different form, its statement/proof can be found in \cite{S} (see Chapter II, Section 5).  Shifted square functions (in various forms) appear naturally in analysis. One such instance refers to the case of the first Calderon commutator, whose approach involves a discrete version of \eqref{Sqf}. This theme is explained in detail in \cite{MS}, Ch.4, pp. 126, while the standard approach for the discretized analogue of \eqref{Sqf} is given at pp. 150.
$\newline$

\noindent\textbf{Step 1. $L^2$ boundedness}
$\newline$

This is a straightforward statement since $\{\phi_{m+j}\}_{j\in\N}$ can be split in at most $C$ disjoint families such that
the functions within each family have pairwise disjoint supports. Thus, from Parseval, we conclude

$$\|S_{l}f\|_2^2=\sum_{j\in\Z}\|(f*\check{\phi}_{j})(x-\frac{l}{2^{j}})\|_2^2$$
$$=\sum_{j\in\Z}\|\hat{f}(\xi)\,\hat{\phi}_{j}(\xi)\,e^{i\,\xi\,\frac{l}{2^{j}}}\|_2^2 \lesssim\|f\|_2^2\;.$$

$\newline$

\noindent\textbf{Step 2. $L^{1}$ to $L^{1,\infty}$ boundedness}
$\newline$

Our intention is to show that there exists $C>0$ absolute constant such that
\beq\label{Sqf1}
\forall\:\:\l>0\:\:\:\:\;\;\:\:\:\:|\{x\,|\,|S_l f(x)|>\l\}|\leq \frac{C}{\l}\,\log (|l|+10)\,\|f\|_1\;.
\eeq
As expected, we will make use of the Calderon-Zygmund decomposition of $f$ at level $\l$:

Thus, for simplicity, we consider $M_d f$ - the \textit{dyadic} Hardy-Littlewood maximal function associated to $f$. Then,
for $\l$ as above, we set $$E_{\l}:=\{x\,|\,M_d f(x)> \l\}\:.$$
  From the definition of $E_{\l}$ we know that there exists a unique collection $\I$ of maximal (disjoint) dyadic intervals such that
\beq\label{El}
E_{\l}=\bigcup_{J\in \I} J\:.
\eeq
  Of course, we retain for later that
\beq\label{meas}
|E_{\l}|=\sum_{J\in \I} |J|\leq \frac{1}{\l}\,\sum_{J\in \I} \int_{J} |f|\leq \frac{1}{\l}\,\|f\|_1\:,
\eeq
and that
\beq\label{Linfb}
|f(x)|\leq \l\:\:\:\:\:\:\:\forall\:\:x\in\R\setminus E_{\l}\;.
\eeq
Next, we decompose the function $f$ as follows:
\beq\label{CZ}
f=g+b\,,
\eeq
with
\beq\label{g}
g:=f\,\chi_{E_{\l}}\,+\sum_{J\in\I}\left(\frac{1}{|J|}\,\int_{J}f\right)\,\chi_{J},
\eeq
and
\beq\label{bbad}
b:=\sum_{J\in\I} b_{J}\:\:\:\:\textrm{with}\:\:\:b_{J}:= \left(f-\frac{1}{|J|}\,\int_{J}f\right)\,\chi_{J}\,.
\eeq
Notice that from the above decomposition we have the following important properties
\beq\label{prop}
\eeq
\begin{itemize}
\item $\|g\|_{\infty}\lesssim \l\,;$
\item $\|g\|_{1}\leq \|f\|_1\,;$
\item $\textrm{supp}\,b_J\subseteq J$ for every $J\in\I$;
\item $\int_{\R} b_{J}=0$ for every $J\in\I$;
\item $\|b_J\|_1\lesssim \l\,|J|$ for every $J\in\I$.
\end{itemize}
Now as in the classical theory approach we continue with the remark
\beq\label{levset}
|\{S_{l}f>\l\}|\leq |\{S_{l}g>\frac{\l}{2}\}|\,+\,|\{S_{l}b>\frac{\l}{2}\}|\;.
\eeq
For the first term we just use the $L^2$ boundedness of $S$ and then the $L^{\infty}-$control on $g$:
\beq\label{gg}
\eeq
$$|\{S_{l}g>\frac{\l}{2}\}|\lesssim \frac{1}{\l^2}\,\|S_{l}g\|_{L^2}^2\lesssim\frac{\|g\|_2^2}{\l^2}\lesssim \frac{1}{\l}\|f\|_1\;.$$
In order to treat the second term, we adopt the following convention: given an interval $J$ with center $c(J)$ and $a>0$, we refer to $a\,J$ as the interval with the same center $c(J)$ and of length $a\,|J|$.

With these done, we have
\beq\label{b1}
\eeq
$$|\{x\,|\,S_{l}b(x)>\frac{\l}{2}\}|$$
$$=|\{x\in\bigcup_{J\in\I}100 J\,|\,S_{l}b(x)>\frac{\l}{2}\}|\,+\,
|\{x\in\R\setminus\bigcup_{J\in\I}100 J\,|\,S_{l}b(x)>\frac{\l}{2}\}|=A\,+\,B\,.$$
Now, on the one hand
\beq\label{A}
A\lesssim \sum_{J\in\I} |J|\leq \frac{\|f\|_1}{\l}\;,
\eeq
while on the other hand
\beq\label{B}
B\lesssim \frac{1}{\l}\,\int_{\R\setminus\bigcup_{J\in\I}100 J}\,S_{l}b\lesssim\frac{1}{\l}\sum_{J\in\I}\int_{\R\setminus100J}
S_{l} b_J\;.
\eeq
It will be thus enough to show that
\beq\label{S}
\int_{\R\setminus100J}S_{l} b_J\lesssim \log (|l|+10)\,\l\,|J|\;.
\eeq
In fact we will show the stronger statement
\beq\label{Ss}
I_{b_J}:=\int_{\R\setminus100J}\sum_{j\in\Z}
|(b_{J}*\check{\phi}_{j})(x-\frac{l}{2^{j}})|\,dx\lesssim \log(|l|+10)\,\l\,|J|\;.
\eeq
For this, we first use the mean zero condition for $b_J$ and write
\beq\label{meanz}
\eeq
$$I_{b_J}^{l,j}(x):=|(b_{J}*\check{\phi}_{j})(x-\frac{l}{2^{j}})|\leq$$
$$\left|\int_{\R}\left(\check{\phi}_{j}(x-c(J)-s-\frac{l}{2^{j}})-\check{\phi}_{j}(x-c(J)-\frac{l}{2^{j}})\right)\,b_{\tilde{J}}(s)\,ds\right|\,,$$
where here we denoted with $c(J)$ the center of $J$, $\tilde{J}=J-c(J)$ and $b_{\tilde{J}}(s):=b_{J}(s+c(J))$.

As a consequence, we have
\beq\label{Ibj}
\begin{array}{cl}
I_{b_J}^{l,j}(x)\lesssim\\\\
\int_{\R}\min\left\{\max_{|u|\leq 2^j\,|J|}\,\frac{2^{2j} |J|}{|2^{j}(x-c(J))-l-u|^2+1},\,\frac{2^{j} }{|2^{j}(x-c(J)-s-\frac{l}{2^{j}})|^2+1}\right\}\,|b_{\tilde{J}}(s)|\,ds\,.
\end{array}
\eeq
Then, implementing \eqref{Ibj} into the LHS of \eqref{Ss}, we have
\beq\label{ctrl}
\eeq
$$|I_{b_J}|\lesssim \l\,|J|\,\sum_{{j\in\Z}\atop{2^{j} |J|\leq 1}} 2^{j} |J|\,\int_{\R\setminus100J}
\frac{2^{j} }{|2^{j}(x-c(J))-l|^2+1}\,dx$$
$$+\,\sum_{{j\in\Z}\atop{2^{j} |J|>1}} \int_{\R\setminus100J}\int_{\tilde{J}}
\frac{2^{j} }{|2^{j}(x-c(J)-s-\frac{l}{2^{j}\,|J|}|J|)|^2+1}\,|b_{\tilde{J}}(s)|\,dx\,ds$$
$$\lesssim  \l\,|J|\,\log (|l|+10 )\,.$$
This finishes Step 2.
$\newline$

\noindent\textbf{Step 3. $L^{q}$ boundedness, $1<q<\infty$}
$\newline$

The fact that \eqref{Sqfp} holds for $1<q\leq 2$ is just a consequence of real interpolation between the $q=2$ case (Step 1)
and $q=1$ case (Step 2).

The other part of the range, \textit{i.e.} $2<q<\infty$, follows from duality. This is again a standard argument, an application of Hincin inequality (probabilistic method). For completeness, we give here the details:

 For $\{\o_{j,l}\}_{{j\in\Z}}$ a sequence of i.i.d. random variables on $[0,1]$ with each $\o_{j}$ taking the values $\{\pm 1\}$ with equal probability (the probability measure here can be taken to be the Lebesgue measure on $[0,1]$) we define for any function $h\in L^1_{loc}(\R)$ the linear operator
\beq\label{Ll}
\L_{l,\o(t)}h(x):=\sum_{j\in\Z}\o_{j}(t)\,(h*\check{\phi}_{j})(x-\frac{l}{2^{j}})\:.
\eeq
  Then, for any $f\in L^q(\R)$, we have
$$\|S_{l}f\|_q^q\approx\int_{[0,1]}\int_{\R}|\L_{l,\o(t)}f(x)|^q\,dt\,dx$$
Taking now any function $g\in L^{q'}(\R)$ with $\|g\|_{q'}=1$ we define
$$V_{l,t}(f,g):=\int_{\R}g(-x)\,\L_{l,\o(t)}f(x)\,dx\,,$$
and notice that
$$V_{l,t}(f,g)=\int_{\R}f(x)\,\L_{l,\o(t)}g(-x)\,dx\,,$$
from which we deduce
\beq\label{defV}
|V_{l,t}(f,g)|\lesssim \|f\|_q\,\|\L_{l,\o(t)}g\|_{q'}\;.
\eeq
Finally, we notice that following exactly the same procedure as that presented at the Step 1 and 2, we have that
\beq\label{Some}
\|\L_{l,\o(t)}h\|_{2}\lesssim \|h\|_{2}\:\:\:\textrm{and}\:\:\:\|\L_{l,\o(t)}h\|_{1,\infty}\lesssim \log(|l|+10)\,\|h\|_{1}\;.
\eeq
Thus from interpolation we get that for any $1< q'\leq2$ we have
\beq\label{Ssq}
\|\L_{l,\o(t)}g\|_{q'}\lesssim_{q'} \log(|l|+10)^{\frac{2}{q'}-1} \|g\|_{q'}\;,
\eeq
with the constant in the above inequality \textit{independent} of $t$.

Now it only remains to observe that we can choose a function $g_t(\cdot)\in L^{q'}(\R)$ such
that $(\int_{[0,1]\times\R}|g_t(x)|^{q'}\,dx\,dt)^{\frac{1}{q'}}=1$ and
\beq\label{Ssqfin}
\|S_{l}f\|_q\approx\int_{[0,1]}\int_{\R}g_t(x)\L_{l,\o(t)}f(x)\,dt\,dx\lesssim_{q'} \log(|l|+10)^{\frac{2}{q'}-1} \|f\|_{q} \|g_t\|_{L_{t,x}^{q'}}\;.
\eeq
\end{proof}

Proposition \ref{Contf} follows now trivially from combining Lemmas \ref{canc}, \ref{max} and \ref{sqf}.
$\newline$
$\newline$
\noindent\textbf{Proof of Proposition \ref{Contgh}.}
$\newline$

The proof of this result follows from Fefferman-Stein's inequality (\cite{FS}) and the result below
$\newline$

\begin{l1}\label{dual} The following holds:
\beq\label{claim}
\sum_{p_0=2^m}^{2^{m+1}}|(g*\check{\phi}_{j,p_0})(x)\,(h*\check{\phi}_{j,p_0})(x)|^{2}\lesssim \|g\|_{\infty}^2\, M(h*\check{\phi}_{\g,m,j})^2(x)\:,
\eeq
where here we set
\beq\label{gmj}
\phi_{\g,m,j}(\eta):=\phi(\frac{\g'(2^{-j})}{2^{j+m}}\,\eta)\;.
\eeq
\end{l1}

\begin{proof}
Notice first that, recalling \eqref{L2infcont}, we have
\beq\label{l2}
\sum_{p_0=2^m}^{2^{m+1}}|(g*\check{\phi}_{j,p_0})(x)|^{2}\lesssim \|g\|_{\infty}^2:.
\eeq
We also know that
\beq\label{MAx}
\begin{array}{cl}
\sup_{p_0\approx2^m}|h*\check{\phi}_{j,p_0}|(x)\\\\
=\sup_{p_0\approx2^m}\left|\int_{\R} (h*\check{\phi}_{\g,m,j})(x-y)\,\frac{2^j}{\g'(2^{-j})}\,\check{\phi}(\frac{2^j}{\g'(2^{-j})}\,y)\,
e^{i\frac{2^j}{\g'(2^j)}\,p_0\,y}\,dy\right|\\\\
\lesssim M(h*\check{\phi}_{\g,m,j})(x)\;.
\end{array}
\eeq
Thus combining \eqref{MAx} and \eqref{l2} we deduce that \eqref{claim} holds.
\end{proof}
$\newline$

\section{Controlled bounds for $\ll_{m}^{+}(f,g,h)$ on the edge $(AB)$.}
$\newline$

 As in the previous section, we only intend here to get $\log$-type upper bounds for our form $\ll_{m}^{+}(f,g,h)$ when restricted to edge $(AB)$ (see Figure 1). It is now simple to notice that taking advantage on the symmetry of our form, by just switching the role of $g$ and $h$ in Theorem \ref{tame}, one gets the desired result:
$\newline$

\begin{t1}\label{tameAB} Let $1<p<\infty$. Then the following holds:
\beq\label{AB}
\underline{\textrm{Edge}\:\:(AB)}\,:\:\:\:\:\:\:\:\:\:\:\:|\ll_{m}^{+}(f,g,h)|\lesssim_{p}  (1+m^{\frac{2}{p'}-1})\,\|f\|_{p}\,\|g\|_{p'}\,\|h\|_{\infty}\,\;.
\eeq
\end{t1}
$\newline$

\section{The proof of our Main Theorem}
$\newline$

It now remains to prove our Theorem \ref{estimabac}. We will do this, by applying real interpolation methods between the tame bounds obtained in Theorems \ref{tame} (and respectively \ref{tameAB}) and the corresponding bounds obtained in \cite{L} for the $L^2$ case.

\subsection{Estimates for $\ll_{m}^{+}(f,g,h)$ at points $B_1$ and $C_1$}

In this subsection we only quote the previous results on the author addressing the (local) $L^2$ case.

As in the previous section, it is useful to notice that because of the symmetry between the role of $g$ and $h$ we only need to refer to one of the mid points of the segments $(AB)$ and $(AC)$ respectively.

However, the estimate for  $\ll_{m}(f,g,h)$ for the point $C_1$ is precisely the content of the author's paper (\cite{L}). Thus, based on Propositions 1 and 2 in \cite{L} together with the symmetric role played by the points $C_1$ and $B_1$ we conclude the following

\begin{t1}\label{estimLblc} \textit{The following estimates hold}:
\beq\label{pointLc}
\underline{\textrm{Point}\:\:C_1(\frac{1}{2},\,\frac{1}{2},\,0)}\,:\:\:\:\:\:\:\:\:\:\:\:|\ll_{m}(f,g,h)|\lesssim 2^{-\frac{m}{16}}\,\|f\|_{2}\,\|g\|_{2}\,\|h\|_{\infty}\;,
\eeq
and
\beq\label{pointLb}
\underline{\textrm{Point}\:\:B_1(\frac{1}{2},\,0,\,\frac{1}{2})}\,:\:\:\:\:\:\:\:\:\:\:\:|\ll_{m}(f,g,h)|\lesssim 2^{-\frac{m}{16}}\,\|f\|_{2}\,\|g\|_{\infty}\,\|h\|_{2}\;.
\eeq
\end{t1}

\subsection{Interpolation: Proof of the main result}

The remaining task for proving our Main Theorem relies on the standard multi-linear interpolation theory.

Indeed, we first interpolate
\begin{itemize}
\item along the edge $(AB)$ between the mid point $C_1$ and a generic point $P\in(AB)$, that is between the bounds obtained
in Theorem \ref{estimLblc} and Theorem \ref{tameAB}.
\item along the edge $(AC)$ between the mid point $B_1$ and a generic point $P\in(AC)$, that is between the bounds obtained
in Theorem \ref{estimLblc} and Theorem \ref{tame}.
\end{itemize}
Through this process we get precisely the content of Theorem \ref{estimabac}. Applying one more time real interpolation between the newly better bounds just obtained on $(AB)$ and $(AC)$, we get the final global result

\begin{t1}\label{final} \textit{Let $P(\frac{1}{p},\,\frac{1}{q},\,\frac{1}{r'})$ be any point in the plane determined by the points $A,\,B,\,C$ such that $P\in \Omega=\overline{\textrm{int}\,\triangle\,ABC}\setminus \{[BC]\cup \{A\}\}$. Then there exists $\a(P)>0$ such that
for any $m\in\N$ the following holds}:
\beq\label{fm}
\underline{\textrm{Point}\:\:P(\frac{1}{p},\,\frac{1}{q},\,\frac{1}{r'})}\,:\:\:\:\:\:\:\:\:\:\:\:|\ll_{m}^{+}(f,g,h)|\lesssim_{P} 2^{-\a(P)\,m}\,\|f\|_{p}\,\|g\|_{q}\,\|h\|_{r'}\;.
\eeq
\end{t1}

Now as mentioned at the end of Section 2, the same reasonings applied to the form
\beq\label{frm}
\ll_m^{-}(f,g,h)=\sum_{j\in\Z\setminus\N}\ll_{j,m}(f,g,h)\,.
\eeq
will give the corresponding analogue of the bound \eqref{fm}.

Finally, recalling that
\beq\label{frm}
\ll(f,g,h)=\sum_{m\in\N}\ll_m(f,g,h)\,
\eeq
and using the geometric summation in the parameter $m$ as resulting from bound \eqref{fm} we deduce that our Main Theorem is true.

\end{document}